\documentclass[12pt,a4paper]{article}

\usepackage{amsmath,amssymb,amsfonts}
\usepackage{epsfig}

\usepackage[latin1]{inputenc}
\usepackage[T1]{fontenc}
\usepackage{ae}
\usepackage[english]{babel}

\usepackage{fancyhdr}
\usepackage{lastpage}
\usepackage{geometry}
\numberwithin{equation}{section}

\newcounter{item}

\renewcommand{\theitem}{\arabic{section}.\arabic{item}}

\newcommand{\cc}{\setcounter{equation}{0}}

\newenvironment{theo}[1]{
\setcounter{item}{\value{equation}}
\addtocounter{equation}{1}
\refstepcounter{item}
\par\addvspace{\bigskipamount}
\indent {\bf \theitem.\hspace{1em}THEOREM#1.} \sl }
{\par\addvspace{\bigskipamount}
}

\newenvironment{pf}{
\par\addvspace{-\smallskipamount}
\indent {\bf Proof.}$\,\ $ } { {$ \, $}
\par\addvspace{\bigskipamount}
}

\newenvironment{lem}{
\setcounter{item}{\value{equation}} \addtocounter{equation}{1}
\refstepcounter{item}
\par\addvspace{\bigskipamount}
\indent {\bf \theitem.\hspace{1em}LEMMA.\,\ }  \sl }
{\par\addvspace{\bigskipamount}
}

\newenvironment{cor}{
\setcounter{item}{\value{equation}} \addtocounter{equation}{1}
\refstepcounter{item}
\par\addvspace{\bigskipamount}
\indent {\bf \theitem.\hspace{1em}COROLLARY.\,\ }  \sl }
{\par\addvspace{\bigskipamount}
}

\newenvironment{ex}{
\setcounter{item}{\value{equation}}\addtocounter{equation}{1}
\refstepcounter{item}
\par\addvspace{\bigskipamount}
\indent {\bf \theitem.\hspace{1em}Example.$\, \ $} }
{\par\addvspace{\bigskipamount} }

\newenvironment{rem}{
\setcounter{item}{\value{equation}}\addtocounter{equation}{1}
\refstepcounter{item}
\par\addvspace{\bigskipamount}
\indent {\bf \theitem.\hspace{1em}{\it Remark.$\, \ $}} }
{\par\addvspace{\bigskipamount}
}

\newenvironment{defi}{
\setcounter{item}{\value{equation}} \addtocounter{equation}{1}
\refstepcounter{item}
\par\addvspace{\bigskipamount}
\indent {\bf \theitem.\hspace{1em}{\it Definition.\,\ }}  \sl }
{\par\addvspace{\bigskipamount} }

\newenvironment{prf}[1]{
\par\addvspace{\bigskipamount}
\indent {\bf #1.$\, \ $} } {\par\addvspace{\bigskipamount} }

\newenvironment{theor}{\bf
\setcounter{item}{\value{equation}}\addtocounter{equation}{1}
\refstepcounter{item}
\par\addvspace{\bigskipamount}
\indent \theitem\hspace{1em}\ignorespaces } {\unskip \ \ \ }

\newcounter{minutes}
\divide\time by 60
\newcounter{hours}
\multiply\time by 60
\addtocounter{minutes}{-\time}

\begin{document}

\begin{center}
{\Large \bf On normal families of quasiregular mappings}
\end{center}
\medskip

\begin{center}
{\large \bf Shamil Makhmutov and Matti Vuorinen}
\end{center}
\bigskip

\medskip

{{\bf Abstract.}
 We discuss the value distribution of quasimeromorphic
mappings in ${\mathbb R}^n$ with given behavior in a neighborhood of
an essential singularity.}

{{\bf 2000 Mathematics Subject Classification.} Primary 30D45.
Secondary 30C65.}

\bigskip

\section{INTRODUCTION}{}

\smallskip

In this paper we study some value distribution properties of
$K$-quasimeromorphic mappings in ${\mathbb R}^n$, $n \geq 2$, with
an essential singularity at infinity and given asymptotic behavior.
We also consider the mutual arrangement of $a$-points of
$K$-quasimeromorphic mappings under certain conditions. In
particular, we consider Yosida type quasimeromorphic mappings (cf.
\cite{12}). This class of mappings was mentioned in \cite[Lemma
1]{1}, \cite{6}, but as far as we know, no further study was carried
out. We prove that $n$-periodic $K$-quasimeromorphic mappings
\cite{5} are Yosida mappings.
\par

\begin{theor}{PRELIMINARIES.} \end{theor}
In general, we follow the notation and terminology of \cite{11}.
The Euclidean distance is denoted by  $|x-y|$ for $x, y \in
{\mathbb R}^n$, and the one-point compactification of ${\mathbb R}^n$ is
denoted by $ {\overline{\mathbb R}}^n\,$ and, as usual, identified
with the Riemann sphere $S^n(\frac{1}{2}e_{n+1},\frac{1}{2})$, \cite[p. 4]{11}.
The chordal distance $q(a , b)$ for $a, b \in
{\overline{\mathbb R}}^n$ is defined by
$$
q(a, b)= \frac{|a-b|}{\sqrt{1+|a|^2}\sqrt{1+|b|^2}} , \quad
\mathrm{if} \quad  a, \, b \, \in {\mathbb R}^n,  \quad \mathrm{ and} \, \,
q(a, \infty)= \frac{1}{\sqrt{1+|a|^2}} \, .
$$
Further let $B_p(a, r)=\{ x \in {\mathbb R}^n : |x-a| <r|a|^{2-p}
\}$, where $a \in {\mathbb R}^n \setminus \{0\}$, $p \geq 1$ and
$r>0$, and $B_2(a,r)=B(a,r)$.
\par

\begin{theor}{NORMAL FAMILIES AND YOSIDA MAPPINGS.} \end{theor}
Here normality of a family of mappings in a domain
$D \subseteq {\mathbb R}^n$ means uniform  continuity
of the family in the spherical metric $q(.,.)$ on
compact subsets of $D$.
\par
We consider the case of families of $K$-quasimeromorphic mappings
in ${\mathbb R}^n$ with essential singularity at infinity only.
\par
\begin{def} \label{def1} \end{def} A family ${\mathcal F} = \{f\}$ of
$K$-quasimeromorphic mappings of a domain
$D \subseteq {\mathbb R}^n$
is called normal on $D$ if every sequence $\{f_n\} \in {\mathcal F}$
has a subsequence that converges uniformly on compact subsets of $D$
with respect to the spherical metric $q$.

\medskip

The criteria for normality of a family of $K$-quasimeromorphic
mappings, proved by R. Miniowitz \cite{6}, can be formulated for $D
\subseteq {\mathbb R}^n$ as follows

\begin{theo}{}\label{thm1} A family ${\mathcal F} =\{f\}$ of
$K$-quasimeromorphic mappings in a domain $D$ in ${\mathbb R}^n$ is
normal if and only if for every compact subset $G$ of $D$ there
exists a number $M_G$ such that
\begin{equation}\label{cond1}
q(f(x_1), f(x_2)) \leq M_G |x_1 -x_2|^{\alpha}
\end{equation}
for each $x_1 \in G$, $x_2 \in D$ and $f \in {\mathcal F}$ where
$\alpha = K^{{1}/({1-n})}$.
\end{theo}

\par
Let
$$
Q_f(x) = \limsup\limits_{|h| \rightarrow 0}
\frac{q(f(x+h), f(x))}{|h|^{\alpha}} \,  .
$$
Then we can rewrite condition (\ref{cond1}) in the form
\begin{equation}\label{cond2}
\sup\limits_{f \in {\mathcal F}} \sup\limits_{x \in G} Q_f(x) \leq
M_G < \infty \ .
\end{equation}
\par
A.Ostrowski \cite{8} applied the technique of normal families of
meromorphic functions for investigation of the value distribution
of meromorphic functions $f$ forming normal families $\{f(2^n
z)\}$ on ${\mathbb C} \setminus \{ 0 \}$ (see also Montel's
monograph \cite{7}). J.Heinonen - J.Rossi \cite{3} and P.J\"arvi
\cite{4} generalized these results on the case of exceptional
$K$-quasimeromorphic mappings in the sense of Julia and described
these mappings in terms of cercles de remplissage.
\par
K.Yosida \cite{12} considered a classification of meromorphic
functions on the complex plane ${\mathbb C}$ based on normality of a
family of functions of the form $\{f(z+a)\}$, $a \in {\mathbb C}$, on
${\mathbb C}$. The paper \cite{13} of L. Zalcman is a review of
normal families of analytic functions.
\par
We define a class of Yosida $K$-quasimeromorphic mappings in
${\mathbb R}^n$ and describe the distribution of $a$-points for these
mappings.
\par

\begin{defi}{} A $K$-quasimeromorphic mapping
$f : {\mathbb R}^n \rightarrow {\overline{\mathbb R}}^n$ is a Yosida
$K$-quasimeromorphic mapping in ${\mathbb R}^n$ if the family of
mappings $\{f(x+a)\}$, $a \in {\mathbb R}^n$, is normal in ${\mathbb
R}^n$.
\end{defi}
\par
The next result follows from Theorem \ref{thm1} with condition
(\ref{cond2}).
\begin{theo}{}\label{thm2} A $K$-quasimeromorphic mapping $f :{\mathbb
R}^n \rightarrow {\overline{\mathbb R}}^n$
is a Yosida mapping in ${\mathbb R}^n$ if and only if
$$
  \sup\limits_{x \in {\mathbb R}^n} Q_f(x) < \infty \ .
$$
\end{theo}

\begin{defi}{}
Let $n \geq 2$ and $p>1$. A $K$-quasimeromorphic mapping $f:
{\mathbb R}^n \rightarrow {\overline{\mathbb R}}^n$
is called $p$-Yosida mapping if for every sequence of
points $\{a_n\}$ in ${\mathbb R}^n$ with $\lim\limits_{n \to \infty}
|a_n|= \infty$, the family of mappings $\{f(a_n+|a_n|^{2-p}x)\}$
is normal in ${\mathbb R}^n$.
\end{defi}

\par
Let $p > 1$ and $y_a=a+|a|^{2-p}x$. Then
$$
Q_f(a)=\lim\limits_{y_a \to a}
\frac{q(f(y_a),f(a))}{|y_a-a|^{\alpha}}=\lim\limits_{|x| \to 0}
\frac{q(f(a+|a|^{2-p}x),f(a))}{| |a|^{2-p}x|^{\alpha}} \ .
$$
Setting  $f_a(x)=f(a+|a|^{2-p}x)$, for every fixed
$a \in {\mathbb R}^n \setminus \{0\}$ we have
$$
Q_{f_a}(0) = \lim\limits_{|x| \to 0}
\frac{q(f_a(x),f_a(0))}{|x|^{\alpha}}
$$
$$
= \lim\limits_{|x| \to 0}
\frac{q(f(a+|a|^{2-p}x),f(a))|a|^{(2-p)\alpha}}{| |a|^{(2-p)}
x|^{\alpha}} =  |a|^{(2-p)\alpha} Q_f(a)  \ .
$$
By Theorem \ref{thm1} with condition (\ref{cond2}), we obtain that
the condition
\begin{equation}\label{cond3}
Q_f= \limsup\limits_{|a| \to \infty} |a|^{(2-p)\alpha} Q_f(a) <
\infty
\end{equation}
is equivalent to the normality of the family of mappings
$\{f(a+|a|^{2-p}x)\}$ in ${\mathbb R}^n$.
\par
Thus we have obtained more general result.

\begin{theo}{} Let $p>1$ and $f$ be a $K$-quasimeromorphic
mapping in ${\mathbb R}^n$. Then $f$ is a $p$-Yosida mapping; that
is, $f$ satisfies {\rm (\ref{cond3})}, if and only if the family of
mappings $\{f(a+|a|^{2-p}x)\}$, $a \in {\mathbb R}^n \setminus
\{0\}$, is normal in ${\mathbb R}^n$.
\end{theo}

Let $\{a_n\}$ be an arbitrary sequence of points in ${\mathbb R}^n$ with
$\lim\limits_{n \to \infty} |a_n| = \infty$, and let $f$ be a
$p$-Yosida mapping in ${\mathbb R}^n$, $p>1$. Without loss of
generality we may suppose that the sequence of mappings
$\{f(a_n+|a_n|^{2-p}x)\}$ converges uniformly to a mapping $g(x)$
on compact subsets of ${\mathbb R}^n$.
\par
Obviously, for every Yosida $K$-quasimeromorphic mapping $f$,
non-constant limit mappings $g$ of converging sequences of the
form $\{f(x+a_n)\}$ are Yosida mappings.

\begin{theo}{} Let $p>1$ and let $f$ be a $p$-Yosida mapping in ${\mathbb
R}^n\,.$ Then non-constant limit mappings of sequences of
the form $\{f(a_n+|a_n|^{2-p}x)\}$, where $\lim\limits_{n \to
\infty} |a_n| = \infty$, are Yosida $K$-quasimeromorphic mappings.
\end{theo}

\begin{pf} Let $e_a=\frac{a}{|a|}$, $a \neq 0$, and $f$ be a
$p$-Yosida $K$-quasimeromorphic mapping in ${\mathbb R}^n$, $p>1$.
Then for every sequence of points $\{a_n\}$, $a_n \in {\mathbb R}^n
\setminus \{0\}$, $\lim\limits_{n \to \infty} |a_n| = \infty$, the
sequence of mappings $\{f(a_n+|a_n|^{2-p}x)\}$ converges uniformly
to a mapping $g$ on compact subsets of ${\mathbb R}^n$, and for
every fixed $x \in {\mathbb R}^n$
$$
Q_{g}(x)= \lim\limits_{y \to x} \frac{q(g(y),
g(x))}{|y-x|^{\alpha}}
$$
$$
\leq \lim\limits_{y \to x} \limsup\limits_{n \to \infty}
\frac{q(f_n(y), f_n(x))}{|y-x|^{\alpha}}
$$
$$
=\lim\limits_{y \to x} \limsup\limits_{n \to \infty}
\frac{q(f(a_n+|a_n|^{2-p}y), f(a_n+|a_n|^{2-p}x))}
{|y-x|^{\alpha}|a_n|^{(2-p)\alpha}} |a_n|^{(2-p)\alpha}
$$
$$
\leq \limsup\limits_{n \to \infty} Q_f(a_n+|a_n|^{2-p}x)
 |a_n|^{(2-p)\alpha}
$$
$$
 =\limsup\limits_{n \to \infty}
Q_f(a_n+|a_n|^{2-p}x) |a_n+|a_n|^{2-p}x|^{(2-p)\alpha}
\frac{|a_n|^{(2-p)\alpha}}{|a_n+|a_n|^{2-p}x|^{(2-p)\alpha}}
$$
$$
\leq Q_f \limsup\limits_{n \to \infty}
|e_{a_n}+|a_n|^{1-p}x|^{(p-2)\alpha}   \, .
$$
Since $p>1$ and $x$ belongs to a compact subset $G \subset {\mathbb
R}^n$, the last limit is bounded. The boundedness of $Q_{g}$
yields that $g$ is a Yosida $K$-quasimeromorphic mapping. {$ \, $}
\end{pf}

\begin{rem}{} {\rm If $p=1$ then we have the class of exceptional
$K$-quasimeromorphic mappings in the sense of Julia (see \cite{3} and
\cite{4}). This kind of mappings form a normal family of mappings $\{f(r_nx)\}$
on ${\mathbb R}^n \setminus \{0\}$ for any sequence of positive numbers $\{r_n\}$,
$\lim\limits_{n \to \infty} r_n= \infty$.}
\end{rem}
\par
\cc
\section{ Cercles de remplissage and $M$-sequences}{}\label{sec2}
\bigskip

S. Rickman \cite{10} proved that for every integer $n \geq 2$ and
each $K \geq 1$ a non-constant $K$-quasimeromorphic mapping $f$ in
${\mathbb R}^n$ can have at most $l(n, K)$ exceptional values where
$l(n, K) \geq 2$.
\par
\begin{defi}{}\label{def4}
Let $p \geq 1$. A sequence $\{x_m\}$ in ${\mathbb R}^n$,
$\lim\limits_{m \to \infty} |x_m| = \infty$, is called a
$M_p$-sequence for a non-constant $K$-quasimeromorphic mapping $f$
in ${\mathbb R}^n$ if for each subsequence $\{x_{m_k}\}$ and each
$\delta >0$ the mapping $f$ takes all but possibly $l=l(n,K)$ values
from  ${\overline{\mathbb R}}^n$ infinitely often in
$\bigcup\limits^{\infty}_{k=1} B_p(x_{m_k}, \delta)$.
\end{defi}
\par
For further investigation we should determine a weighted distance
$d_p(X, Y)$ between two sequences of points $X=\{x_m\}$ and
$Y=\{y_m\}$ in ${\mathbb R}^n \setminus\{0\}$, where $\lim\limits_{m
\to \infty} |x_m| = \infty$ and $\lim\limits_{m \to \infty} |y_m| =
\infty$. Let
$$
D_p(X, Y)= \inf_{k,m}\dfrac{|x_m-y_k|}{|x_m|^{2-p}}\, ,
\quad
D_p(Y, X)= \inf_{k,m}\dfrac{|y_k-x_m|}{|y_k|^{2-p}} \, .
$$
\begin{lem}{}
Let $p \geq 1$. For any two sequences of points $X=\{x_m\}$ and
$Y=\{y_m\}$ in ${\mathbb R}^n \setminus \{0\}$ with $\lim\limits_{m
\to \infty} |x_m| = \infty$ and $\lim\limits_{m \to \infty} |y_m| =
\infty$, either $D_p(X,Y)$ and $D_p(Y,X)$ both are zero or both are
non-zero.
\end{lem}
\par
\begin{pf} Indeed, if $D_p(X,Y)=0$, that is, there are
subsequences of points $\{x_{m_j}\}$ and $\{y_{k_j}\}$ that
satisfy the conditions
$$
|x_{m_j}-y_{k_j}|<\varepsilon_j  |x_{m_j}|^{2-p}, \quad  \text{
where} \quad \varepsilon_j \to 0 \, \quad  \text{ as} \quad j \to
0.
$$
Then
$$
\dfrac{|x_{m_j}-y_{k_j}|}{|y_{k_j}|^{2-p}} =
\dfrac{|x_{m_j}-y_{k_j}|}{|x_{m_j}|^{2-p}}
\dfrac{|x_{m_j}|^{2-p}}{|y_{k_j}|^{2-p}} \leq
   \varepsilon_j \dfrac{|x_{m_j}|^{2-p}}{|y_{k_j}|^{2-p}}
\quad .
$$ Since $p \geq 1$ we see that
$$
||x_{m_j}|^{2-p}-|y_{k_j}|^{2-p}| \leq M ||x_{m_j}|-|y_{k_j}||
\leq M |x_{m_j}-y_{k_j}| <  \varepsilon_j |x_{m_j}|^{2-p} \quad .
$$
This inequality implies that
$\left|1-\dfrac{|y_{k_j}|^{2-p}}{|x_{m_j}|^{2-p}} \right|< \varepsilon_j$.
According to the last estimate $|x_{m_j}| \sim |y_{k_j}|$ eventually, i.e.
$$
||x_{m_j}|-|y_{k_j}|| < \varepsilon_j
|x_{m_j}|^{2-p} \leq \hat\varepsilon_j |y_{k_j}|^{2-p}  \, .
$$
Thus we conclude that $D_p(X,Y)=0$ implies  $D_p(Y,X)=0$ as well.
{$ \, $}
\end{pf}
For  two sequences of
points $X$ and $Y$ in ${\mathbb R}^n\setminus \{0\}$ set
$$
d_p(X, Y)= \min \{D_p(X,Y), D_p(Y,X)\} \, .
$$
\par
\begin{theo}{}\label{thm4}
Let $p>1$ and $f$ be a $K$-quasimeromorphic mapping in ${\mathbb
R}^n$, $n \geq 2$. A sequence of points $\{x_m\}$ in ${\mathbb
R}^n$, $\lim\limits_{m \to \infty} |x_m| = \infty$, is a
$M_p$-sequence for $f$ if and only if for each $c \in
{\overline{\mathbb R}}^n$ there is a sequence of points $\{x'_m\}$
in ${\mathbb R}^n$ such that $d_p(\{x_m\},\{x'_m\}) = 0$ and
$\lim\limits_{m \to \infty} q(f(x'_m),c) = 0$.
\end{theo}
\par
\begin{pf} Necessity. If $\{x_m\}$ is a
$M_p$-sequence for $f$ then for each $k>0$ there is a finite number
of points $x_m$ such that the image of the balls
$B_p(x_m,\frac{1}{k})$ is spherically bounded away from the point $c
\in {\overline{\mathbb R}}^n$ by the distance $\frac{1}{k}$. That
is, for each $k>0$ there is a number $m_k$ such that for each
$m>m_k$ there is a point $x'_m \in B_p(x_m,\frac{1}{k})$ with
$q(f(x'_m), c) < \frac{1}{k}$. For $m \leq m_1$ we take $x'_m=x_m$,
and for $m_k < m \leq m_{k+1}$ let $x'_m$ be as above. The sequence
of points $\{x'_m\}$ is then "$p$-closed" to $\{x_m\}$, i.e. ,
$\lim\limits_{m \to \infty}\dfrac{|x_m-x'_m|}{|x_m|^{2-p}}=0$, and
$q(f(x'_m),c) \rightarrow 0$ as $m \rightarrow \infty$.
\par
Sufficiency. We take sequences of points $\{x'_k\}$ and $\{x_k\}$
in ${\mathbb R}^n$ such that
$$d_p(x'_k, x_k)\rightarrow 0 \, \,  \mathrm{ as} \, \,k
\rightarrow \infty \,\, \mathrm{ but} \, \,
q_p(f(x'_k), f(x_k))\nrightarrow 0 \,\, \mathrm{as} \,\,k
\rightarrow \infty$$
at the same time. Consider a family of
mappings ${\mathcal F} = \{f(x_k+|x_k|^{2-p}x) \}$. By our assumption
the family ${\mathcal F}$ is not equicontinuous in the neighborhood of
the origin. Therefore, for each $r>0$ the family ${\mathcal F}$ is not
normal in the ball $B_r=\{x \in {\mathbb R}^n , \, |x|\leq r\}$ (see
\cite{11}). By Theorem 4 \cite{6}, the family $\mathcal F$ takes
infinitely often all but possibly $l(n, K)$ values from
${\overline{\mathbb R}}^n$ in the ball $B_r$. Therefore the
mapping $f$ takes infinitely often all but possibly $l(n, K)$
values from ${\overline{\mathbb R}}^n$ in the union of balls
$\bigcup\limits^{\infty}_{k=1} B_p(x_k, r)$. Since the mentioned
argument is valid for all $r>0$ and each subsequence $\{x_k\}$ we
can conclude that $\{x_m\}$ is a $M_p$-sequence for $f$. {$ \, $}
\end{pf}
\par
\begin{cor}{} If $\{x_m\}$ in ${\mathbb R}^n$ with
$\lim\limits_{m \to \infty} |x_m| = \infty$, is a $M_p$-sequence for
a mapping $f$ in ${\mathbb R}^n$ then a sequence $\{y_m\}$ in
${\mathbb R}^n$ satisfying the condition of closeness
$\lim\limits_{m \to \infty}\dfrac{|x_m-y_m|}{|y_m|^{2-p}}=0$ is also
$M_p$-sequence for $f$.
\end{cor}
\par
\begin{defi}{}\label{def5}
Let $p \geq 1$. A sequence $\{x_m\}$ in ${\mathbb R}^n$ with
$\lim\limits_{m \to \infty} |x_m| = \infty$, is called
$\mu_p$-sequence for a $K$-quasimeromorphic mapping $f$ if there
exist monotone decreasing sequences of positive numbers $\{L_m\}$,
$\lim\limits_{m \to \infty} L_m=0$, and $\{r_m\}$, $\lim\limits_{m
\to \infty} r_m=0$, such that the image of each ball $B_p(x_m,
r_m)$, $m=1, 2, \dots$, covers the sphere $S^n$ except possibly
$l(n, K)$ sets $E_j$, $j= 1, \dots, l(n,K)$, whose spherical
diameters do not exceed $L_m$.
\end{defi}

\begin{rem}{} {\rm The concepts of $M_p$-sequence and
$\mu_p$-sequence for meromorphic functions on the complex plane
${\mathbb C}$ were introduced by V.Gavrilov in \cite{2}.}
\end{rem}
\par

\begin{lem}{}\label{lemma2}
A sequence of points $\{x_m\} \in {\mathbb R}^n$,
$\lim\limits_{m \to \infty} |x_m| = \infty$,  is a $\mu_p$-sequence
for a $K$-quasimeromorphic mapping $f$ in ${\mathbb R}^n$, $p>1$, if and only
if for each $r>0$ there exist sets $E_j(r,m) \in S^n$,
$j=1, \dots, l(n,K)$, such that the spherical diameter of the sets do
not exceed $r$, and $N(r)$ is such that the mapping $f$
takes in every ball $B_p(x_m,r)=\{x \in {\mathbb R}^n, \, |x-x_m|<r
|x_m|^{2-p} \}$, $m \geq N$,  all values from ${\overline{\mathbb R}}^n$
except possibly the sets $E_j(r, m)$, $j=1, \dots, l(n,K)$.
\end{lem}
\par
\begin{pf} The proof of necessity follows from the definition of
$\mu_p$-sequence.
\par
Sufficiency.  Let $\{x_m\}$ be a sequence that satisfies the
conditions of the Lemma. Let $r=\frac{1}{l}$ where $l=1, 2, \dots$.
Consider the corresponding sets $E_j(\frac{1}{l}, m)$, $l=1, 2,
\dots$; $m=1, 2, \dots$, and $j=1, \dots, l(n,K)$. Then we obtain
the sequence of numbers $\{N(\frac{1}{l})\}$ such that
$$
N\left(\frac{1}{1}\right) < N\left(\frac{1}{2}\right) < \dots
< N\left(\frac{1}{l}\right) < \dots \, \, .
$$
Now we assume that $r_m=\frac{1}{l}$ and $L_m=\frac{1}{l}$ for all
$m$ from the interval $(N(\frac{1}{l}), N(\frac{1}{l+1}))$. For $m
\leq N(\frac{1}{1})$, we assume $r_m=1$ and $L_m=1$, i.e. the
diameter of the sphere $S^n$. Choosing the sequences $\{r_m\}$ and
$\{L_m\}$ according to described scheme, we obtain that $\{x_m\}$
satisfies the definition of $\mu_p$-sequences of points for
$K$-quasimeromorphic mappings. {$ \, $}
\end{pf}

\begin{lem}{}\label{lemma3} Let $p > 1$. Then $f$ is a $p$-Yosida
$K$-quasimeromorphic mapping in ${\mathbb R}^n$ if and only if $f$
does not possess a $M_p$-sequence.
\end{lem}
\par
\begin{pf} Necessity. Let $f$ be the $p$-Yosida mapping and
$\{x_m\}$ be the $M_p$-sequence of points for $f$. Consider the
sequence of mappings $\{f(x_m+|x_m|^{2-p}x)\}$ in $\{|x|\leq r\}$.
Since the sequence is normal in $\{|x|\leq r\}$, without loss of
generality we suppose that it converges on $\{|x|\leq r\}$. There are
two possible cases: either the limit mapping of the sequence is a
constant (including infinity) or a non-constant mapping.
\par
In the first case, one can take $r$ sufficiently small and $N$
sufficiently big such that the mappings $f(x_m+|x_m|^{2-p}x)$
assume the values closed to the limit value in $\{|x|\leq r\}$ for
all $m \geq N$. In other words we can formulate this fact as
following: $f$ assumes the values closed to the limit value in
each ball $B_p(x_m,r))$, $m \geq N$, and therefore, the images
$f(B_p(x_m,r))$ do not cover $S^n$. Hence, $\{x_m\}$ cannot be a
$M_p$-sequence for $f$.
\par
Now we suppose that $g$ is a non-constant $K$-quasimeromorphic limit
mapping of the sequence $\{f(x_m+|x_m|^{2-p}x)\}$ in $\{|x|\leq
r\}$. We choose $r$ such that $q(g(x),g(0)) \leq 1/4$ in $\{|x|\leq
r\}$. According to the Hurwitz theorem for quasimeromorphic mappings
\cite[Lemma 2]{6}, \cite{9}, there is a number $N$ such that
$q(f(x_m+|x_m|^{2-p}x),g(x)) \leq 1/2$ in $\{|x|\leq r\}$ for $m
\geq N$. The last implies that $\{x_m\}$ cannot be a $M_p$-sequence
for $f$.
\par
Sufficiency. If $f$ does not possess $M_p$-sequences of points for
$f$ then  for every sequence $\{x_k\}$ with $\lim\limits_{k \to
\infty}|x_k|= \infty$, and every $\delta >0$ we always can find a
subsequence $\{x_m\}$ such that $f$ assumes in
$\bigcup\limits^{\infty}_{m=1} B_p(x_m, \delta)$ all values from
${\overline{\mathbb R}}^n$ except possibly at least $l(n, K)+1$
values. In other words, the family $\{f(x_m+|x_m|^{2-p}x)\}$ assumes
in $\{|x|\leq \delta\}$ all values from ${\overline{\mathbb R}}^n$
except possibly at least $l(n, K)+1$ values. According to
\cite[Theorem 5]{6}, the family $\{f(x_k+|x_k|^{2-p}x)\}$ is normal
in $\{|x|\leq \delta\}$. Since $\{x_k\}$ is an arbitrary sequence,
$f$ is a $p$-Yosida quasimeromorphic mapping. {$ \, $}
\end{pf}

Next lemma shows equivalence of the notions of $\mu_p$-sequences
and $M_p$-sequences for spatial $K$-quasimeromorphic mappings.
\par
\begin{lem}{} Let $p > 1$. A sequence $\{x_m\}$ in ${\mathbb
R}^n$, $\lim\limits_{m \to \infty} |x_m| = \infty$, is a
$\mu_p$-sequences of points for a $K$-quasimeromorphic mapping $f$
in ${\mathbb R}^n$ if and only if $\{x_m\}$ is a $M_p$-sequence for
$f$.
\end{lem}
\begin{pf} It follows from definitions \ref{def4} and \ref{def5} that each
$\mu_p$-sequences of points for a $K$-quasimeromorphic mapping $f$
is a $M_p$-sequence for $f$.
\par
To prove the converse statement we suppose that $\{x_m\}$ is a
$M_p$-sequence for $f$ and it is not a $\mu_p$-sequence. Then there
is a positive $\varepsilon_0$ for which Lemma \ref{lemma2} is not
satisfied, i.e. there is a subsequence $\{x_k\}$ in ${\mathbb R}^n$
such that for each $k$ the set of values in ${\overline{\mathbb
R}}^n$ not assumed by $f$ in the ball $B_p(x_k, \varepsilon_0)$
cannot be contained in $l$ sets whose spherical diameters do not
exceed $r_0$. It implies that $f$ omits $l+1$ distinct values $c_{1,
k}$, $c_{2, k}$, \dots , $c_{l+1,k}$ from ${\overline{\mathbb R}}^n$
in each ball $B_p(x_k, \varepsilon_0)$ and
\begin{equation}\label{cond4}
q(c_{j, k}, c_{i, k}) \geq \frac{r_0}{2} \quad , \quad i \neq j
\end{equation}
Consider the family of mappings $\{f(x_k+|x_k|^{2-p}x)\}$ in the
ball $B_0=\{|x|\leq \varepsilon_0\}$. Since every mapping
$f_k(x)=f(x_k+|x_k|^{2-p}x)$ does not take $l+1$ values satisfying
(\ref{cond4}) in the ball $B_0$, then by Theorem 5 \cite{6}, the
family $\{f_k(x)\}$ is normal in $B_0$. But according to our
assumption, $\{x_m\}$ is a $M_p$-sequence for $f$ and,
consequently by Lemma \ref{lemma3}, the family of mappings cannot
be normal at the origin. That contradicts our assumption. {$ \, $}
\end{pf}

Now we consider the main result of the paper: mutual arrangement
of $a$-points for $p$-Yosida $K$-quasimeromorphic mappings in
${\mathbb R}^n$.
\par
\begin{theo}{}\label{thm5} Let $f$ be a $K$-quasimeromorphic mapping $f$
in ${\mathbb R}^n$. Let $A_1, \dots, A_{l+2}$, where $l=l(n, K)$,
be distinct points in ${\overline{\mathbb R}}^n$ and
$f^{-1}(A_j)=\{a_{jk} \vert \, k \in {\mathbb N} \}$, $j=1, 2,
\dots, l+2$. Then $f$ is a $p$-Yosida $K$-quasimeromorphic mapping
in ${\mathbb R}^n$, $p>1$ if and only if
\begin{equation}\label{cond5}
\inf\limits_{i \neq j, k, m} \left\{
\frac{|a_{jk}-a_{im}|}{|a_{jk}|^{2-p}} \, : \, \, a_{jk} \neq 0
\right\} >0 \, .
\end{equation}
\end{theo}
\par
\begin{pf} Necessity. Suppose that the limit value of (\ref{cond5}) for
sequences $\{a_{1 k}\}$ and $\{a_{2 m}\}$ is equal to 0. Then by
Theorem \ref{thm4}, both sequences are $M_p$-sequences for $f$
and, consequently, $f$ is not a $p$-Yosida mapping.
\par
Sufficiency. Suppose that $f$ is not a $p$-Yosida
$K$-quasimeromorphic mapping in ${\mathbb R}^n$. Then there is a
$\mu_p$-sequence of points $\{x_m\}$ for the mapping $f$. By the
definition, there exists a number $N$ such that from $N$ onward
every ball $B_p(x_m, \varepsilon_m)$,  $\lim\limits_{m \to \infty}
\varepsilon_m =0$, contains the roots of at least two equations
$f(x)=A_1$ and $f(x)=A_2$, where $A_1, A_2 \in {\overline{\mathbb
R}}^n$. We choose accordingly those $A_1$-points and $A_2$-points,
that belong to $B_p(x_m, \varepsilon_m)$ and denote them as
$\{a_{1 m}\}$ and $\{a_{2 m}\}$, resp. Then these points satisfy
the conditions
$$
 |a_{1 m}-x_m|< \varepsilon_m |x_m|^{2-p} \, ,
$$
$$
|a_{2 m}-x_m|< \varepsilon_m |x_m|^{2-p} \,  ,
$$
and hence
$$
\lim\limits_{m \to \infty} \frac{|a_{1 m}-a_{2 m}|}{|a_{1
m}|^{2-p}} = 0  \, .
$$
The last contradicts the conditions of the theorem.
{$ \, $}
\end{pf}
\par
The result similar to Theorem \ref{thm5} for $p=1$; that is, Julia
exceptional functions, was obtained by P.J\"arvi \cite[Theorem
6]{4}.
\par


\cc
\section{ Yosida quasimeromorphic mappings of the first order}{}


Similar to K.Yosida \cite{12}, we define Yosida $K$-quasimeromorphic
mappings of the first order.
\par
\begin{defi}{} A Yosida $K$-quasimeromorphic  mapping
$f$ in ${\mathbb R}^n$ is called Yosida $K$-quasimeromorphic mapping
of the first order if for every sequence of points $\{x_n\}$ in
${\mathbb R}^n$ with $\lim\limits_{n \to \infty} |x_n|= \infty$,
limit mappings of converging subsequences of $\{f(x+x_n)\}$ on
compact subsets of ${\mathbb R}^n$ are not constants.
\end{defi}
\par
For every set $E \subset {\overline{\mathbb R}}^n$, we denote by $q(E)$
the spherical diameter of $E$.

\begin{theo}{}\label{thm6} Let $f$ be a Yosida $K$-quasimeromorphic
mapping.

{\rm (a)}   For every $0<d<1$ there exists $r>0$ such that
   $$ \sup\limits_{x \in {\mathbb R}^n}  q(f(\overline{B}(x, r)))\leq d \, .$$

{\rm (b)} If for every $r>0$ there exists $\delta >0$, such that
   \begin{equation}\label{cond6}
   \inf\limits_{x \in {\mathbb R}^n}  q(f(\overline{B}(x, r))) \geq \delta \,
   \end{equation}
   then
   \begin{equation}\label{cond7}
   \sup\limits_{a \in \overline{\mathbb R}^n}  n_f (x, r, a) \leq C \, ,
   \end{equation}
   where $C$ depends of $r$ and $f$.

\end{theo}
\par
\begin{pf} Condition (a) follows from the definition of Yosida
$K$-quasimeromorphic mappings and Theorem \ref{thm2}.
\par
To prove condition (b) we assume that there exist sequences of
points $\{x_k\}$ in ${\mathbb R}^n$ with$\lim\limits_{k \to \infty}
|x_k|=\infty$, and $\{a_k\}$ in ${\overline{\mathbb R}}^n$,
$\lim\limits_{k \to \infty} |a_k|=a$, satisfying the condition
$\lim\limits_{k \to \infty} n_f(x_k, r, a_k)=\infty$. Consider a
family $\{f_k(x)\}$, $f_k(x)=f(x+x_k)$, which is normal in ${\mathbb
R}^n$. Without loss of generality, we assume that $\{f_k(x)\}$ is a
converging sequence on compact subsets of ${\mathbb R}^n$, and by
(\ref{cond6}), the limit mapping $F$ is non-constant. Since $F$ is a
discrete mapping we may choose $R$ such that $a \notin F(\partial
B(0, R))$. The sequence $\{f_k\}$ converges uniformly to $F$ in
$\overline{B}(0, R)$ and therefore
$$
f_k(\partial B(0,R)) \cap \{a, a_k\} = \emptyset \, .
$$
Hence, by a property of the topological index [9, p.86]
$$
n_f(x_k, r, a_k) \leq n_{f_k}(0, R, a_k) =n_{f_k}(0, R, a) =n_F(0, R, a) < \infty
$$
for large $k$. This contradicts our assumption.
{$ \, $}
\end{pf}

\begin{cor}{} If $f$ is a Yosida $K$-quasimeromorphic
mapping and for every $r>0$ there exists $\delta >0$ such that
$$
\inf\limits_{x \in {\mathbb R}^n}  q(f(\overline{B}(x, r))) \geq \delta
$$
then $f$ is a Yosida $K$-quasimeromorphic mapping of the first
kind.
\end{cor}

\begin{ex}{} {\rm $n$-periodic quasimeromorphic mappings are
Yosida quasimeromorphic mappings of the first kind. Such mappings
were studied by O. Martio and U. Srebro in \cite{5}. }
\end{ex}
Recall that for a quasimeromorphic mapping $f$ in $B(0, R)$ and
$r<R \leq \infty$ the average counting function $A_f(r)$ is defined as
$$
A_f(r)=\frac{1}{\lambda_n} \int\limits_{\overline{\mathbb R}^n}
\frac{n(0,r,y)}{(1+|y|^2)^n} dm(y) \, ,
$$
where
$$
\lambda_n=\int_{\overline{\mathbb R}^n}
\frac{1}{(1+|y|^2)^n} dm(y) \,
$$
and $m(y)$ is a Lebesgue measure in $\overline{\mathbb R}^n$.

\begin{theo}{}\label{thm7} Let $f$ be a Yosida $K$-quasimeromorphic
mapping in ${\mathbb R}^n$. Then
$$
A_f(r)=O(r^n)  \ \ \  {\text as} \ \ \ r \to \infty \, .
$$
Furthermore, if $f$ is a Yosida $K$-quasimeromorphic mapping of
the first order then
\begin{equation}\label{cond8}
0<\lim\limits_{r \to \infty} \frac{A_f(r)}{r^n} < \infty \, .
\end{equation}
\end{theo}
Let
$$
A_f(x,r)=\frac{1}{\lambda_n} \int\limits_{\overline{\mathbb R}^n}
\frac{n(x,r,y)}{(1+|y|^2)^n} dm(y) \, .
$$

\begin{lem}{}\label{lemma5} If $f$ is a Yosida $K$-quasimeromorphic
mapping of the first order then for every fixed $r>0$
$$
\liminf\limits_{|x| \to \infty} A_f(x,r) > 0 \, .
$$
\end{lem}
\begin{pf}
Assume that $\liminf\limits_{|x| \to \infty} A_f(x,r)=0$. Then there
exists a sequence $\{x_k\}$ in ${\mathbb R}^n$ on which $n(x_k,r,y)
\to 0$ as $k \to \infty$ for almost all $y \in \overline{\mathbb
R}^n$. Therefore $A_f(x_k,r)=A_{f_k}(r) \to A_F(r)
>0$ as $k \to \infty$ since $F \neq const$. This is a contradiction.
{$ \, $}
\end{pf}

\begin{prf}{Proof of Theorem \ref{thm7}}
If $f$ is a Yosida $K$-quasimeromorphic mapping, then by Theorem
\ref{thm6} (a), for every $0<d<1$ there exists $\rho
>0$ such that $q(f(\overline{B}(x,\rho))) \leq d$ for every $x \in
{\mathbb R}^n$ and therefore  $A_f(x,\rho) \leq C < \infty$.
\par
Let us consider the ball $B(0,r)$ where $r$ is sufficiently large,
and Whitney decomposition $\mathcal W$ of $B(0,r)$. That is, $\mathcal W$
is a union of dyadic essentially disjoint closed cubes $Q_k$ of
the diameter $\rho$ and centered at $a_k$. Each cube $Q_k$ can
be inscribed into the ball $B(a_k,\rho/2)$. Then $B(0,r) \subset
\bigcup\limits_{k} B(a_k,\rho/2)$ and therefore
$$
\begin{array}{ll}
     A_f(r)&  = \int\limits_{\overline{\mathbb R}^n} \frac{n(0,r,y)}{(1+|y|^2)^n}\, dm(y)
     = \int\limits_{B(0,r)} \frac{J(x,f)}{(1+|f(x)|^2)^n}\, dm(x) \\ &
     \leq \sum_{k} \int\limits_{B(a_k, \rho/2)} \frac{J(x,f)}{(1+|f(x)|^2)^n}\, dm(x)
     \leq C N(r)
      \, ,
\end{array}
$$
where $N(r)$ is the number of cubes $Q_k$ in the Whitney
decomposition $\mathcal W$. One can see that $N(r) \thickapprox
\dfrac{\text{Volume } B(0, r)}{\text{Volume } B(x, \rho/2)} =
\left(\dfrac{2}{\rho}\right)^n r^n$. Hence, $A_f(r)= O(r^n)$ as $r
\to \infty$.
\par
To prove (\ref{cond8}) we will use a similar method. We take the
cubes $Q_k$ in the Whitney decomposition $\mathcal W$ such that
each ball $B(x, \rho)$ is inscribed into $Q_k$. By Lemma
\ref{lemma5}, $A_f(x, r) \geq \delta
>0$ for every $x \in {\mathbb R}^n$. Then
$$
\begin{array}{ll}
     A_f(r)&  = \int\limits_{\overline{\mathbb R}^n} \frac{n(0,r,y)}{(1+|y|^2)^n}\, dm(y)
     = \int\limits_{B(0,r)} \frac{J(x,f)}{(1+|f(x)|^2)^n}\, dm(x) \\ &
     \geq \sum_{k} \int\limits_{B(a_k, \rho)} \frac{J(x,f)}{(1+|f(x)|^2)^n}\, dm(x)
     \geq \delta N(r)
      \, .
\end{array}
$$
Thus,  $\dfrac{A_f(r)}{r^n} >0$ and (\ref{cond8}) is proved. {$ \,
$}
\end{prf}


\cc
\section{ Further remarks and questions }{}

In this final part of the paper we make some remarks and pose some natural questions.
\par
In \cite{5}, O.~Martio and U.~Srebro (Theorem 7.4)  proved that if
$f : {\mathbb R}^n \rightarrow {\overline{\mathbb R}}^n$ is an
$n$-periodic $K$-quasimeromorphic mapping then
\begin{equation}
N(f, y, \tilde F_f) = \sum\limits_{x \in f^{-1}(y)\cap F_f} i(x,f)
< \infty
\end{equation}
for all $y \in {\overline{\mathbb R}}^n$,
where $\tilde F_f$ is a fundamental set for $f$.
\par
It follows from Theorems \ref{thm5} and \ref{thm6} that mentioned
result can be generalized on Yosida $K$-quasimeromorphic mappings.
In other words, we can say: if $f$ is a Yosida $K$-quasimeromorphic
mapping then for every $r>0$
$$
\sup\limits_{a \in {\mathbb R}^n} N(f, y, B(a,r)) < \infty
$$
for all $y \in {\overline{\mathbb R}}^n$.
\par
It will interesting to prove the following conjectures.

\begin{theor}{Problem.} {\rm $f : {\mathbb R}^n \rightarrow
{\overline{\mathbb R}}^n$ is a Yosida $K$-quasimeromorphic mapping
if and only if there is $r>0$ such that
$$
\sup\limits_{a \in {\mathbb R}^n} N(f, y, B(a,r)) < \infty
$$
for all $y \in {\overline{\mathbb R}}^n$. Furthermore,
$$
\inf_{y \in {\overline{\mathbb R}}^n}
\sup\limits_{a \in {\mathbb R}^n} N(f, y, B(a,r)) > 0
$$
 is the necessary and sufficient condition for $f$ to be
 a Yosida $K$-quasimeromorphic mapping.
}
\end{theor}

In Theorem \ref{thm7}, we gave the estimation of the growth of
$A_f(r)$ for Yosida quasimeromorphic mappings and Yosida
$K$-quasimeromorphic mappings of the first order. Since Yosida
mappings are characterized as mappings with uniform behavior in
${\mathbb R}^n$, it is a natural question to describe Yosida
$K$-quasimeromorphic mappings in terms of uniform growth of
$A_f(r)$ in ${\mathbb R}^n$.

{\bf Acknowledgements.}
Research supported in part by Sultan Qaboos University, project
IG/SCI/DOMS/04/04 and by the Academy of Finland.

\bigskip

\noindent

\textit{Shamil Makhmutov} \textbf{\hfill} E-mail: {\tt makhm@squ.edu.om }
\newline
ADDRESS: {\it Department of Mathematics and Statistics, Sultan Qaboos University, P.O. Box 36, Al Khodh 123, Oman;
Institute of Mathematics, Ufa, 450077, Russia}\\

\medskip
\textit{Matti Vuorinen} \textbf{\hfill} E-mail: {\tt vuorinen@utu.fi}
\newline
ADDRESS: {\it Department of Mathematics, FIN-20014 University of Turku, Finland}


\small
\begin{thebibliography}{99}


\bibitem [1] {1} {\sc   A.~Eremenko}: {\em  Bloch radius, normal families
and quasiregular mappings},
 Proc. Amer. Math. Soc. 128(2) (1999),   557--560.


\bibitem [2] {2} {\sc V.~I.~Gavrilov}: {\em  Behavior of a meromorphic
function in a neighborhood of its essential singularity.} (Russian)
 Izv. Akad. Nauk SSSR Ser. Mat. 30 (4), 1966, 767--788.


\bibitem [3] {3} {\sc J.~Heinonen and J.~Rossi}: {\em Remarks on the value
distribution of quasimeromorphic mappings.}
Complex Variables Theory Appl. 21( 1993), 231--242.

\bibitem [4] {4} {\sc  P.~J\"arvi}: {\em On the behavior of quasiregular
mappings in the neighborhood of an isolated singularity,} Ann.
Acad. Sci. Fenn. Ser. AI Math. 15, 1990, 341--353.

\bibitem [5] {5} {\sc O.~Martio and U.~Srebro}:
{\em Periodic quasimeromorphic mappings. } J.
Anal. Math. 28 (1975), 20--40.



\bibitem [6] {6} {\sc R.~Miniowitz}: {\em Normal families of
quasimeromorphic mappings},  Proc. Amer. Math. Soc. 84 (1982), 35--43.

\bibitem [7] {7} {\sc P.~Montel}: {\em Lecons sur les familles normales de
functions analytiques et leurs applications}, Gauthier-Villars,
Paris, 1927.

\bibitem [8] {8} {\sc A.~Ostrowski}:   {\em \"Uber Folgen analytischer
Funktionen und einige Versch\"arfungen des Picardschen Satzes},
Math. Z. 24 (1925),   215--258.

\bibitem [9] {9} {\sc Yu.~G.~Reshetnyak}:
{\em Space mappings with bounded distortion}, Transl. of Math.
Monographs, Vol. 73, AMS, 1989.

\bibitem [10] {10} {\sc S.~Rickman}: {\em  On the number of omitted values
of entire quasiregular mappings}, J. Anal. Math. 37 (1980),
100--117.





\bibitem [11] {11} {\sc M. Vuorinen}: {\em Conformal Geometry
and Quasiregular Mappings}, Lecture Notes in Math. 1319,
Springer-Verlag, Berlin--New York, 1988.



\bibitem[12]{12}   {\sc K.~Yosida}: {\em On a class $(A)$ of
meromorphic functions},  Proc. Phys.-Math. Soc. Japan  16 (1934),
227 -- 235.

\bibitem[13]{13}   {\sc L. Zalcman}: Normal families: new perspectives.
Bull. Amer. Math. Soc. (N.S.)  35  (1998),  no. 3, 215--230.

\end{thebibliography}

\normalsize
\end{document}